\documentclass[12pt, leqno]{article}
\usepackage{amsmath}
\usepackage{amssymb}
\usepackage{amscd}
\usepackage{amsthm}
\def\overset#1#2{{\mathrel{\mathop {{#2}_{}}\limits^{#1}}}}
\def\underset#1#2{{\mathrel{\mathop {{}_{} {#2}}\limits_{{#1}_{}}}}}
\def\upplim_#1{\underset{#1}{\overline\lim}\;}
\def\lowlim_#1{\underset{#1}{\underline\lim}\;}

\newcommand{\C}{{\mathbf{C}}}

\newcommand{\cL}{{\mathcal{L}}}
\newcommand{\cO}{{\mathcal{O}}}

\newcommand{\im}{{\Im\,}}
\newcommand{\iso}{\cong}

\newcommand{\lam}{{\lambda}}

\newcommand{\nb}{{\bf N.B.}  }

\newcommand{\one}{{\mathbf{1}}}

\newcommand{\re}{{\Re\,}}

\newcommand{\supp}{{\mathrm{Supp}\,}}

\newcommand{\Z}{{\mathbf{Z}}}

\newtheorem{lem}[equation]{\bf Lemma}
\newtheorem{thm}[equation]{\bf Theorem}
\newtheorem{claim}[equation]{\it \quad Claim}

\newtheorem{prob}[equation]{\bf Problem}
\newtheorem{prop}[equation]{\bf Proposition}

\newtheorem{quest}[equation]{\bf Question.}

\numberwithin{equation}{section}
\newenvironment{pf}
{{\it Proof}\hskip10pt} {\hfill{\it q.e.d.}\par\vskip+10pt}
\title{On Oka's Extra-Zero Problem
}
\author{Makoto Abe, Sachiko Hamano and Junjiro Noguchi}
\date{Version 1 (2011/aug./10)}
\begin{document}
\maketitle
\thispagestyle{empty}
\baselineskip=18pt
\parindent12pt
\begin{abstract}
After the solution of Cousin II problem by K. Oka III in 1939,
he thought an {\it extra-zero problem} in 1945 (his posthumous paper)
asking if it is possible to solve an arbitrarily given Cousin II problem
adding some extra-zeros whose support is disjoint from the given one.
By the secondly named author, some special case was
affirmatively confirmed in dimension two
and a counter-example in dimension three or more was given.
The purpose of the present paper is to give a complete solution
of this problem with examples and to discuss some new questions.
\end{abstract}

\section{Introduction.}
After the solution of Cousin II problem by K. Oka \cite{oka} III
he thought the following {\it extra-zero problem} in 1945
(his posthumous paper \cite{oka-psth}, No. 2, p.\ 31 Problem 2;
see \S2):
\smallskip

{\it
{\bf Oka's Extra-Zero Problem.}
``Let $X$ be a domain of holomorphy and let
$D$ be an effective divisor on $X$.
Find an effective divisor $E$ on $X$ such that
their supports have no intersection,
$$
(\supp D) \cap (\supp E)= \emptyset ,
$$ and Cousin II problem
for $D+E$ is solvable on $X$.''
}

Let $L(D)$ denote the line bundle determined by $D$,
let $N(D)=L(D)|_{(\supp D)} \to \supp D$ be the normal
bundle of $D$ over the support $\supp D$ of $D$,
and $\one_X$ denote the trivial line bundle over $X$.
Then Cousin II problem is equivalent to ask if $L(D) \iso \one_X$.
{\it Oka Principle} (\cite{oka}, III)
says that  $L(D) \iso \one_X$ if and only if
the the first Chern class $c_1(L(D))=0$ in the cohomology group
$H^2(X, \Z)$.
Since the problem is trivial for $D$ such that $L(D)\iso \one_X$,
Oka's extra-zero problem makes sense for $D$ with $c_1(L(D))\not=0$.

In \cite{ham07} a counter-example was constructed
in $\dim X \geq 3$, and if $\dim X=2$,
some partial affirmative answer was shown.

The purpose of this paper is to give a complete answer
to Oka's extra-zero problem with examples and
some new questions based on this problem, on which
we would like to put equal emphasis as well (cf.\ \S\S 4 and 5).
The following is our main theorem.

\begin{thm}
\label{main}
Let $D$ be an effective Cartier divisor on a Stein space $X$.
Then Oka's extra-zero problem is solvable if and only if
$c_1(N(D))=0$ in $H^2(\supp D, \Z)$.

In particular, if $\dim X=2$, Oka's extra-zero problem is
always solvable.
\end{thm}

The last statement is due to $H^2(\supp D, \Z)=0$,
since $\dim \supp D=1$.

{\nb}  K. Oka \cite{oka-psth} almost proved Theorem \ref{main}
(see Theorem \ref{oka}). Referring to Oka's Theorem \ref{oka},
one may say that
Theorem \ref{main} is an infinitesimalization of the
topological condition from a neighborhood of $D$ to $D$ itself.
This is not difficult now by many well-established results.

The divisor $E$ in Oka's extra-zero problem is called an
{\it extra-zero} of $D$. 
By definition $L(E)=L(-D)$. Thus the problem is equivalent
to find a holomorphic section $\sigma \in \Gamma(X, L(-D))$
such that
\begin{equation}
\label{sigma}
\supp (\sigma) \cap \supp D=\emptyset.
\end{equation}
Here we consider only $\sigma$ whose zero set is
nowhere dense in $X$ and hence defines a divisor $(\sigma)$
on $X$.
From this viewpoint it is interesting to see
\begin{prop}
\label{opp}
Let the notation be as in Theorem \ref{main}.
Then Oka's extra-zero problem is solvable if and only if
there exists a section $\tau \in \Gamma(X, L(D))$
with nowhere dense zero set and
\begin{equation}
\label{tau}
\supp (\tau ) \cap \supp D=\emptyset.
\end{equation}
\end{prop}

\nb For $\tau$ in \eqref{tau} we required that the zero set of
$\tau$ is nowhere dense in $X$. This is, however, not a restriction.
For if $\tau$ vanishes constantly on an irreducible component
$X'$ of $X$, then we take a section $\tau' \in \Gamma(X, L(D))$
such that $\tau'|_{X'} \not\equiv 0$ and
$\tau' \equiv 0$ on every irreducible component of $X$ other than
$X'$. Then $\{\tau+\tau'=0\} \subset \{\tau=0\}$ as sets
and $\{\tau+\tau'\}|_{X'} \not\equiv 0$.
In this way we may modify $\tau$ so that its zero set is nowhere dense
in $X$.  This is the same for $\sigma$ in \eqref{sigma}.
\smallskip

{\it Acknowledgment.}
After the counter-example constructed by \cite{ham07}
which is a reducible divisor, Professor T. Ueda
asked if there is an irreducible counter-example;
his question forms a part of the motivation
of the present paper.
Professor S. Takayama gave an interesting example of \S4.
Professor T. Tsuboi kindly answered to the thirdly named author
a number of questions on the triangulation of
complex analytic subsets.
The authors are very grateful to all of them.

\section{Oka's notes.}
Here we summerize in short the contents of the
posthumous paper \cite{oka-psth}.
We should first notice that it is dated 28 February 1945
before Oka's Coherence Theorem (\cite{oka} VII).

Roughly speaking, he developed the following study.
\begin{enumerate}
\item
He wished to reformulate Cousin II problem
by relaxing the conclusion so that it is solvable
on every domain of holomorphy.
\item
He recalled the Oka Principle for Cousin II problem on
a domain of holomorphy, and reduced the essential
key-part of the problem to the following:

\quad{\it Let $\bar\Omega \Subset \C^n$ be a bounded closed domain with a
holomorphically convex neighborhood.
Let $D$ be divisor on a neighborhood $\bar{\Omega}$.
Then the Cousin II problem for $D$ is solvable in a
neighborhood of $\bar\Omega$
iff $c_1(L(D))=0$ in a neighborhood
of $\bar\Omega$.\footnote{Here his term is ``balayable'' used in Oka
\cite{oka} III; the meaning is that the given Cousin II distribution
is continuously deformable to a zero-free continuous Cousin II distribution.
The Cousin II problem on a domain $X$ of holomorphy
is solvable iff $D$ is ``balayable'' on $X$.
}}
\item
He then posed the Extra-Zero Problem as Problem 2.
Let $\Omega$ and $D$ (effective) be as in the above item.
Then he asks to find an effective divisor $E$ in a
neighborhood of $\bar\Omega$
such that $\supp D \cap \supp E=\emptyset$  and Cousin II problem
for $D+E$ is solvable in a neighborhood of $\bar\Omega$.
\item
He proved a result as Theorem 8 which is stated as follows:

\quad Theorem. {\it The extra-zero problem is solvable for $D$ in
a neighborhood of $\bar{\Omega}$ if and only if there is a
neighborhood $V$ of $\bar{\Omega}$ with
$c_1(L(D)|_V)=0$.}
\item
After confirming the above topological obstruction for
the extra-zero problem, he proved that
there always exists an effective divisor $F$ in a neighborhood
of $\bar{\Omega}$ such that Cousin II problem for $D+F$
is solvable.
Furthermore he proved that there are at most $n+1$ holomorphic
functions $f_j, 1 \leq j \leq n+1$, in a neighborhood of $\bar{\Omega}$
such that in a neighborhood $W$ of every point of $D \cap \bar{\Omega}$
one of zeros of $f_j$ is exactly $D \cap W$.
\end{enumerate}

Taking account of the above items (ii) and (iv),
we may assume that he obtained or at least recognized
the following statement.

\begin{thm} {\rm (Oka \cite{oka-psth}) }
\label{oka}
Let $\Omega \subset \C^n$ be a domain of holomorphy,
and let $D$ be an effective divisor on $\Omega$.
Then the extra-zero problem for $D$ is solvable if and only if
there is a neighborhood $V$ of $D$ satisfying
$c_1(L(D)|_V)=0$ in $V$.
\end{thm}

K. Oka wrote that it strongly attracts his interest
from a number of viewpoints to decide if
this Extra-Zero Problem is always solvable or
there is a counter-example, and the problem would have
a wide influence in future.\footnote{
He did not give an explicit problem here.}

It is now necessary to know what is the most general
form of his statement (Theorem \ref{oka}), and it is
Theorem \ref{main}.

\section{Proofs.}
\hskip12pt {\bf (a) Proof of Theorem \ref{main}.}
Suppose first that Oka's extra-zero problem is solvable.
Let $E$ be an extra-zero of $D$, and let
$\sigma \in \Gamma(X, L(E))$ with $(\sigma)=E$.
Set
$$
Y=\supp D.
$$
Then the restriction $\sigma|_{U}$ to $U=X \setminus \supp E$
has no zero over the neighborhood $U$ of $Y$.
Therefore $L(-D)|_U=L(E)|_U \iso \one_U$, and then
$N(D)\iso \one_{Y}$, so that $c_1(N(D))=0$.

Conversely, assume that $c_1(N(D))=0$.
Note that $c_1(N(D)) \in H^2(Y, \Z)$ is a restriction of
$c_1(L(D)) \in H^2(X, \Z)$.
By Mihalache \cite{mih96} there is a Stein neighborhood
$V$ of $Y$ for which there is a strong deformation retract
$V \to Y$.
Therefore we have
$$
H^2(V, \Z) \iso H^2(Y, \Z).
$$
It follows that $c_1(L(D))|_V=0$.\footnote{
The existence of such a neighborhood $V$ for a chosen element
$c_1(L(D))$ is sufficient for our
argument here.
For that purpose, it suffices to know the two facts:
(1) There is a system of neighborhoods of $Y$ which admit deformation
retracts to $Y$ (Whitney-Bruhat, Comment. Math. Helv. 1959):
(2) There is a system of neighborhoods of $Y$ which are Stein
 (Siu, Invent. Math. 1976).
}
Since $V$ is a Stein space, we have that $L(D)|_V=\one_V$,
and hence $L(-D)|_V=\one_V$.
Thus there is a section $\sigma \in \Gamma(V, L(-D))$
without zero on $V$.
By the Fundamental Theorem of Oka-Cartan
(Oka \cite{oka} I--II, VII--VIII; Grauert-Remmert \cite{gr})
the restriction $\sigma|_{Y}$ extends to a holomorphic
section $\tilde{\sigma} \in \Gamma(X, L(-D))$
with nowhere dense zero set.
Thus the divisor $(\tilde{\sigma})$ gives rise to an extra-zero
of $D$.
\smallskip

{\bf (b) Proof of Proposition \ref{opp}.}
We keep the notation used in the above (a).
Suppose that Oka's extra-zero problem is solvable.
Then the above $\sigma \in \Gamma(X, L(E))$ has no zero
on $Y$.
Therefore,
$N(D)=L(D)|_{Y}=L(-E)|_{Y}\iso \one_{Y}$.
By the Fundamental Theorem of Oka-Cartan
$\sigma^{-1}|_{Y}$ holomorphically extends to
a section $\tau \in \Gamma(X, L(D))$ with nowhere dense zero set.
By definition $\supp (\tau) \cap Y= \emptyset$.

Suppose the existence of $\tau \in \Gamma(X, L(D))$
with nowhere dense zero set
such that $\supp (\tau) \cap Y= \emptyset$.
Then the same argument implies the existence
of $\sigma \in \Gamma(X, L(-D))$ with nowhere dense zero set
such that
$\supp (\sigma) \cap Y= \emptyset$, and hence
$(\sigma)$ is an extra-zero of $D$.

\section{Examples.}

\hskip12pt{\bf (a)} The first solvable non-trivial example for
Oka's extra-zero problem was given by \cite{ham07} Theorem 1.
Using a similar idea we give another example.
Let $X=(\C^*)^2$
with $\C^*=\C\setminus \{0\}$.
Then the torus $T=S^1 \times S^1 \subset X$ gives the generator
of $H_2(X, \Z) \iso H^2(X, \Z)$. K. Stein \cite{st41}, \S4 computed the divisor
$D$ on $X$ corresponding to $T$, i.e., $c_1(L(D))=T$.
Let $(z,w)\in X$ be the natural coordinates.
Then the analytic hypersurface  given by
\begin{equation}
\label{stdiv}
D^+: \qquad w=z^i=e^{i \log z}
\end{equation}
has the first Chern class $T$.
Stein \cite{st41} also obtained an analytic function
$F^+(z, w)$ that defines $D$:
\begin{equation}
\label{stfun}
F^+(z,w)=\exp\left(\frac{(\log z)^2}{4\pi} +\frac{\log z}{1-i}\right)
\prod_{\nu=0}^\infty \left(1-\frac{w}{e^{i\log z+2\nu\pi}}\right)
\times
\prod_{\mu=1}^\infty \left(1-\frac{1}{we^{-i\log z+2\mu\pi}}\right),
\end{equation}
where we take a branch $\log 1=0$.
Then $\langle c_1(L(D^+)), T \rangle =1$, and so Cousin II problem for $D^+$
is not solvable.
Let $\cL_z$ denote the analytic continuation as the variable $z$
runs over the unit circle in the anti-clockwise direction.
Then $\cL_z \log z= \log z +2\pi i$, and
$$
\cL_z F^+(z,w)= w F^+(z,w), \qquad \cL_w F^+(z,w)=F^+(z,w).
$$

Set
\begin{align*}
D^- : \qquad w &=z^{-i}=e^{-i \log z},\\
F^-(z,w) &= F^+\left( \frac{1}{z}, w \right).
\end{align*}
Then $L(D^+ + D^-)\iso \one_X$, however $D^+ \cap D^- \not=\emptyset$.
We have by \eqref{stfun} that
\begin{align}
\label{stfun1}
F^-(z,w) &=\exp\left(\frac{(\log z)^2}{4\pi} -\frac{\log z}{1-i}\right)
\prod_{\nu=0}^\infty \left(1-\frac{w}{e^{-i\log z+2\nu\pi}}\right)
\times
\prod_{\mu=1}^\infty \left(1-\frac{1}{we^{i\log z+2\mu\pi}}\right),\\
\nonumber
& \cL_z F^-(z,w) = \frac{1}{w} F^-(z,w), \qquad \cL_w F^-(z,w)=F^-(z,w).
\end{align}

By Theorem \ref{main} there is an extra-zero $E$ of $D^+$,
but it is unknown what is $E$.
Therefore it is very interesting to ask
\begin{quest}
Find an analytic expression of $E$.
\end{quest}

On the other hand we may give an example for Proposition \ref{opp}.
Let $\lam \in \C$ such that the real part
$\Re \lam \not\in 2 \pi \Z$ and set
$$
D^+_\lam : \qquad w=e^\lam z^{i}.
$$
Then $D^+_\lam \cap D^+= \emptyset$ , $L(D^+_\lam)=L(D^+)$,
and $D^+_\lam$ is the zero of the analytic function
$$
F^+_\lam (z,w)=F^+\left( z, e^{- \lam} w\right).
$$

Set $
\Omega=\left\{\xi \in \C; |\re\, \xi| < {\pi},
|\im \, \xi| < {\pi}
\right\}$.  Then, it is interesting to observe that the holomorphic mapping
\begin{equation}
\Phi : \zeta \in \C \to (e^\zeta, e^{i\zeta})\in (\C^*)^2=X
\end{equation}
is into-biholomorphic; this describes precisely why $D^+$ is
``balayable'' in a neighborhood of $D^+$ (see \S2 (ii) and its footnote).
\smallskip

{\bf (b) Examples for Theorem \ref{main} with
$c_1(N(D))\not=0$.}

{\bf (1) (Reducible divisor)} 
A counter example in $\dim X \geq 3$ is given in \cite{ham07}
in a domain of $\C^n$ ($n \geq 3$).
Using a similar idea, we give another counter example of a divisor
on $(\C^*)^3$ for which Oka's extra-zero problem has
no solution.

Now we let $X=(\C^*)^2 \times \C^*=(\C^*)^3$
with projection $p:X \to (\C^*)^2$.
Let $D^+ \subset (\C^*)^2$ be as in the above (a), and
set
\begin{align}
\label{red}
D_1 &= D^+ \times \C^*, \qquad D_2=(\C^*)^2\times \{1\},\\
\nonumber
D &= D_1+D_2.
\end{align}
Since $L(D_2)\iso \one_X$, $L(D)\iso L(D_1) \iso p^* L(D^+)$.
Therefore $N(D)|_{D_2}\iso L(D^+) \not\iso \one_{D_2}$
with $D_2 \iso (\C^*)^2$,
so that $N(D) \not\iso \one_D$.
One sees that $D$ has no extra-zero on $X$.

{\bf (2) (Irreducible divisor)}
The above example of $D$ is reducible, and we like to
have an irreducible analytic hypersurface that has
no extra-zero. We are going to modify the example of (1).

Let $\Z[i]=\Z+ i\Z$ be the lattice of Gaussian integers and put
$$
\C \to \C/\Z[i]\: \iso \: \C^*
\: \overset{\lam_0}{\to} \: \C^*/\Z=E \: \iso \: \C/\Z[i] .
$$
Then $E$ is an elliptic curve with complex multiplication
$a \in E \to ia \in E$. Set
$$
\iota: (a,b) \in E^2 \to (ia, b) \in E^2
$$
and let $\Delta \subset E^2$ be the diagonal divisor.
Set
\begin{align*}
D_1 &=\iota^*\Delta, \\
\lam_1 & =\lambda_0 \times \lambda_0: E^2 \to E^2, \\
\hat{D}_1 &= \lambda^* D_1 \subset (\C^*)^2.
\end{align*}

Note that the example of Stein \cite{st41}, \S4
($w=z^i$ in $({\C^*})^2$) is
a connected component $\hat{D}'_1$ of $\hat{D}_1$.
It follows that the Chern class
$$
c_1(L(\hat{D}_1)) \not= 0 \quad \mbox{in } H^2((\C^*)^2, \Z).
$$
(This is equivalent to the non-solvability of Cousin II for $\hat{D}_1$,
or to the non-triviality of the line bundle
$L(\hat{D}_1)$ over $(\C^*)^2$.)
In fact, letting $T=S^1 \times S^1 \in H_2((\C^*)^2, \Z)$ denote the
generator, we get
\begin{equation}
\label{2cycle}
\langle c_1(L(\hat{D}_1)), T \rangle =
\langle c_1(L(\hat{D}'_1)), T \rangle = 1.
\end{equation}

Now we set
\begin{align*}
\lam_2 : X & =(\C^*)^3 \times \C^* \to E^2 \times E \quad
(\mbox{the quotient map}),\\
D_2 &= D_1 \times E +E^2 \times \{0\},\\
\hat{D}_2 &=\lam_2^* D_2.
\end{align*}
Then $L(\lam_2^* (E^2 \times \{0\}))$ is trivial on $X$ and so
$L(\hat{D}_2)=L(\lam_2^*(D_1 \times E))$,
which is the pull-back of $L(\hat{D}_1)$ over $(\C^*)^2$
by the projection $X \to (\C^*)^2$.
Therefore, $L(\hat{D}_2) \not\iso \one_{X}$,
where $\one_X$ denotes the trivial line bundle over $X$.

Furthermore, we see that the normal bundle
$N(\hat{D}_2)=L(\hat{D}_2)|_{\hat{D}_2} \to \hat{D}_2$ is non-trivial.
For $N(\hat{D}_2)|_{(\C^*)^2 \times\{1\}}\iso L(\hat{D}_1)$.
Therefore we obtain
\begin{lem}
\label{lem1}
Let the notation be as above. Then
$L(\hat{D}_2) \not\iso \one_X$ and
$N(\hat{D}_2) \not\iso \one_{\hat{D}_2}$.
\end{lem}

\nb  This means that Cousin II problem for $\hat{D}_2$ on $X$
is not solvable and there is no extra zero for $\hat{D}_2$.

We would like to deform $\hat{D}_2$ to a smooth irreducible
divisor, but this is not trivial. Thus we are going to deform $D_2$ on $E^3$,
but $D_2$ is not ample. To make it ample, we add the divisor
$\{1\} \times E^2$ to $D_2$ with setting
$$
D_3=D_2+\{1\} \times E^2,
$$
which is then ample, and we put $\hat{D}_3=\lam_2^*D_3$ on $X$.
Since $\lam_2^* L(\{1\} \times E^2)=
L(\lam_2^{-1}\{1\} \times (\C^*)^2) \iso \one_X$,
$$
L(\hat{D}_3) \iso L(\hat{D}_2).
$$
Thus Lemma \ref{lem1} holds for $\hat{D}_3$, too:
\begin{lem}
\label{lem2}
Let the notation be as above. We have that
$L(\hat{D}_3) \not\iso \one_X$ and
$N(\hat{D}_3) \not\iso \one_{\hat{D}_3}$.
\end{lem}

It is well known that $L(3 D_3)$ is very ample.
We take a smooth irreducible hyperplane section $D_4$ by a holomorphic
section of $L(3D_3)$, and set
$$
\hat{D}_4=\lam_2^* D_4.
$$
\begin{prop}
{\rm (Example) }  Let the notation be as above.
Then $\hat{D}_4$ is a smooth irreducible divisor on $X$
such that $L(\hat{D}_4) \not\iso \one_X$ and
$N(\hat{D}_4) \not\iso \one_{\hat{D}_4}$; equivalently,
\begin{align*}
c_1(L(\hat{D}_4)) &\not=0 \quad \mbox{in } H^2(X, \Z),\\
c_1(N(\hat{D}_4)) &\not=0 \quad \mbox{in } H^2(\hat{D}_4, \Z).
\end{align*}
\end{prop}
\begin{pf}
It is clear due to the construction that $\hat{D}_4$ is smooth and
irreducible (or connected).
Now we look at the 2-cycle $T$ in \eqref{2cycle}.
We regard $T=S^1 \times S^1 \times \{1\} \in H_2(X, \Z)$.
Then this cycle $T$ comes from a 2-cycle of $E^3$, which
is again denoted by the same $T \in H_2(E^3, \Z)$.
Then it follows that
\begin{equation}
\label{pair1}
\langle c_1(L({D}_4)), T \rangle = 3,
\end{equation}
so that $c_1(L(\hat{D}_4)) \not=0$.

It remains to show that $c_1(N(\hat{D}_4)) \not= 0$.
By Lefschetz' hyperplane-section theorem the natural morphism
$$
H_2(D_4, \Z) \to H_2(E^3, \Z) \to 0
$$
is surjective, and then there is a 2-cycle
$T' \in H_2(D_4, \Z)$ which is mapped to $T$.
Then $T'$ can be lifted to a 2-cycle in $H_2(\hat{D}_4, \Z)$, denoted
by the same $T'$.
We see by \eqref{pair1} that
$$
\langle c_1(N(\hat{D}_4)), T' \rangle = \pm 3.
$$
Thus $c_1(N(\hat{D}_4)) \not=0$;
this finishes the proof.
\end{pf}

{\bf (3) (Takayama's irreducible example) }
Let $z_j=x_j+iy_j, 1 \leq j \leq n$ be the natural complex
coordinates of $\C^n$ with the standard basis
$e_j, 1 \leq j \leq n$. Then $e_j, ie_j, 1 \leq j \leq n$ form
real basis of $\C^n$ and we define a lattice $\Gamma \subset \C^n$
defined by
$$
\Gamma=\langle e_1, \ldots, e_n, i e_1, \ldots , ie_n \rangle.
$$
We set $A=\C^n/\Gamma$ and a sequence of covering maps,
\begin{equation}
\label{quot}
\C^n\, \overset{\rho}{\to}\, (\C^*)^n\, \overset{\pi}{\to}\, A,
\end{equation}
where $\rho$ is the quotient map by
$\langle ie_1, e_2, \ldots, e_n \rangle$ and $\pi$ is that by
$\langle e_1, ie_2, \ldots, ie_n \rangle$.
We set $X=(\C^*)^n$.

Let $L$ be the line bundle whose Chern class is represented by
$$
\omega= d i\sum_{j=1}^n dz_j \wedge d\bar{z}_j +
i \sum_{j \not= k} dz_j \wedge d\bar{z}_k,
\qquad d \in \Z.
$$
Then $L$ is ample for $d \geq 2$, and very ample if $d \geq 4$.
\begin{claim}
$\pi^* \omega \not= 0$ in $H^2((X, \Z))$;
in particular, the pairing,
$\omega \cdot (ie_1 \wedge e_j) \not= 0$, $j \geq 2$
where $ie_1 \wedge e_j \in H_2(A, \Z)$.
\end{claim}
{\it Proof.}
We consider the two pull-back morphisms
\begin{equation*}
\pi^* : H^q(A, \Z) \to H^q(X , \Z), \qquad q \geq 1.
\end{equation*}
Then $\pi^* dx_1=0$, and $\pi^*dy_k=0, k \geq 2$;
on the other hand, $\pi^* dy_1 \not= 0$, and $\pi^*dx_k=0, k \geq 2$.
It follows that
$$
i dz_j \wedge dz_j=2 dx_j \wedge dy_j =0
\quad (\mbox{mod } dx_1, dy_k, k \geq 2).
$$
Therefore we have
$$
\pi^* i(dz_1 \wedge d\bar{z}_j+dz_j\wedge d\bar{z}_1)=
\pi^*(i dy_1 \wedge d_j + dx_j \wedge (-idy_1))
=-2 \pi^* (dy_1 \wedge dx_j)
$$
for $j \geq 2$.\hfill $\square$

Now we assume $n \geq 3$ and $d \geq 4$.
Then $L$ is very ample.
\begin{prop} {\rm (Example)}
\label{taka}
We  take a smooth irreducible divisor $D \in |L|$
and set $\tilde{D}=\pi^{-1} D \subset X$.
Then the divisor $\tilde{D}$ is smooth irreducible and
has no extra-zero on $X$.
\end{prop}
\begin{pf}
Since $H_1 (D, \Z) \iso H_1(A, \Z)$ (Lefschetz' Theorem),
$\tilde{D}$ is connected.
Again by Lefschetz' Theorem the natural morphism is surjective:
$$
H_2(D, \Z) \to H_2(A, \Z) \to 0.
$$
There is an element $\xi \in H_2(D, \Z)$ which is mapped to
$ie_1 \wedge e_j$ ($j \geq 2$).

Let $\iota: D \hookrightarrow A$ be the inclusion map and
let $\tilde{\iota}: \tilde{D} \to X$ be the lifting.
It follows from \eqref{quot} that there is an element
$\tilde{\xi} \in H_2(\tilde{D}, \Z)$ with
$\tilde{\iota}_*\tilde{\xi}=\xi$.
Note that $c_1(L(\tilde{D}))=\pi^* \omega$.
We have that
\begin{equation*}
c_1(L(\tilde{D})) \cdot {\tilde{\xi}}=
\omega \cdot (ie_1 \wedge e_j) \not=0.
\end{equation*}
Therefore we see that $c_1(L(\tilde{D}))\not= 0$
and that $c_1(N(\tilde{D}))\not=0$;
equivalently, the smooth irreducible divisor
$\tilde{D}$ has no extra-zero on $X$.
\end{pf}

\section{Intersections of analytic cycles.}

We would like to consider what is the intersection theory
of analytic cycles on Stein manifolds.

{\bf (a)}  The prototype of intersection theory is Bezout's Theorem
such that for two cycles $A_1, A_2$ on a variety
$$
\deg (A_1 \cdot A_2)=\deg A_1 \cdot \deg A_2.
$$
Cornalba-Shiffman \cite{cs} however gave a counter example
of analytic curves $C_1$ and $C_2$ of $\C^2$
such that the ``orders'' of $C_j$ are zero, but the
``order'' of $C_1 \cdot C_2$ is infinite.

From the viewpoint of Oka's extra-zero problem,
however, we should have by Theorem \ref{main}
$$
C_1 \cdot C_2 =0.
$$
More in general, because of Theorem \ref{main}
there is no global intersections
of ``divisors'' and ``curves'' on Stein manifolds
(nor of divisors and analytic cycles on $\C^n$ ($n \geq 2$)).
In fact, let $D$ be a divisor on a Stein manifold $X$ and
let $C$ be an analytic curve in $X$.  Then the intersection
$$
  D \cdot C
= \langle c_1(L(D)), C \rangle
= \langle c_1(L(D)|_C), C \rangle
= \langle c_1(\one_C), C \rangle
= 0 .
$$

This suggests that some topological structure must be involved
in the possible intersection theory on a Stein manifold.
Therefore it is interesting to propose

\begin{quest}
What is the global intersection theory of analytic cycles
on Stein manifolds?
\end{quest}

We may consider at least three kind of intersections
on a Stein manifold $X$:
\begin{enumerate}
\item
Intersections of zeros (divisor) of holomorphic functions $X$:
\item
Intersections of of hypersurfaces of $X$ with trivial normal bundles:
\item
Intersections of hypersurfaces of $X$.
\end{enumerate}

\nb It is noticed that the normal bundle of the zeros (divisor)
of a holomorphic function is trivial.

\begin{quest}
What are the differences of these intersections?
In particular, characterize the complete intersections of hypersurfaces
with trivial normal bundles on a Stein manifold.
\end{quest}

\nb If $X$ is affine algebraic, there are intersections in
algebraic category and there is a difference even in a simplest
case as follows. Let $X \subset \C^2$ be an affine elliptic curve
with a point at infinity, and let $a \in X$ be a point,
which is an algebraic divisor.
There is no regular rational function on $X$ with exact zero $a$,
but there exists such a holomorphic function on $X$.

\begin{quest}
Let $X$ be Stein and algebraic.
Let $D$ be an effective algebraic divisor on $X$ with
$c_1(L(D))=0$ (resp.\ $c_1(N(D))=0$).
Does there exist a holomorphic function $f \in \cO(X)$
with zero divisor $D$ (resp.\ locally in a neighborhood of $\supp D$)
such that the order of $f$ at the infinity is at most one.
\end{quest}

{\bf (b)}
We set $X=(\C^*)^3$ and would like to discuss the global
intersections of analytic cycles on $X$.
As observed in (a), there is no
intersection between analytic curves and divisors on $X$.
therefore we may restrict ourselves to deal with
the intersections of divisors on $X$.
The first homology group of $X$ is
$$
H_1 (X, \Z) \iso \Z^3,
$$
which is generated by $e_1=S^1 \times \{1\}^2$,
$e_2=\{1\}\times S^1 \times \{1\}$ and $e_3=\{1\}^2 \times S^1$.
Then their products generate the higher homology groups and
in particular,
$$
H^2(X, \Z) \iso H_2(X, \Z) \iso \Z^3.
$$
\smallskip

{\bf (c)}  {\it Stein's example from the viewpoint of the value distribution
theory.}
Let $f: \zeta \in \C \to (e^\zeta, e^{i\zeta}) \in (\C^*)^2=X$
be the example \eqref{stdiv} due to Stein in \S4.
Then $f$ is algebraically non-degenerate; that is, there is no
proper algebraic subset $Y \subset X$ with $f(\C) \subset Y$.
In fact, let $P(z,w) (\not=0)$ be any non-zero
polynomial in $(z,w) \in X$.
We write
$$
P(z,w)=\sum_{j,k} c_{j,k}z^j w^k.
$$
Suppose that $f(\C) \subset \{P=0\}$. Then
$$
\sum_{j,k} c_{j,k} e^{(j+ik)\zeta} \equiv 0.
$$
This is absurd, since
$e^{(j+ik)\zeta}$ are linearly independent over $\C$.

According to the main result of Noguchi-Winkelmann-Yamanoi \cite{nwy02},
\cite{nwy08}, and Corvaja-Noguchi \cite{cn}, the intersection
set $f(\C) \cap D$ is infinite for an arbitrary algebraic divisor $D$
on $X$. For an extra-zero $E$ of $D^+=f(\C)$ we have
$$
f(\C) \cap E= \emptyset.
$$
\begin{prob}
Let $g: \C \to X$ be an analytically non-degenerate entire curve.
Then, is $g(\C) \cap A \not= \emptyset$ for an
arbitrary analytic divisor $A$ of $X$?
Moreover, is $g(\C) \cap A$ an infinite set?

Here it is natural to generalize $X$ to a semi-abelian variety.
\end{prob}

\rightline{Division of Mathematical and Information Sciences}
\rightline{Faculty of Integrated Arts and Sciences}
\rightline{Hiroshima University}
\rightline{1-7-1 Kagamiyama}
\rightline{Higashi-Hiroshima, 739-8521 Japan}
\rightline{e-mail: abem@hiroshima-u.ac.jp}
\medskip

\rightline{Department of Mathematics}
\rightline{Faculty of Human Development and Culture}
\rightline{Fukushima University}
\rightline{Kanayagawa, Fukushima, 960-1296 Japan}
\rightline{e-mail: hamano@educ.fukushima-u.ac.jp}
\medskip

\rightline{Graduate School of Mathematical Sciences}
\rightline{The University of Tokyo}
\rightline{Komaba, Meguro, Tokyo, 153-8914 Japan}
\rightline{e-mail: noguchi@ms.u-tokyo.ac.jp}

\end{document}